\documentclass[cleveref,a4paper,numberwithinsect]{lipics-v2021}
\pdfoutput=1

\usepackage{cite} 
\usepackage{amsmath}
\usepackage{amssymb}  
\usepackage{cleveref}
\usepackage{enumerate}
\usepackage{tikz}
\tikzstyle{every picture} = [>=latex]
\usetikzlibrary{calc,arrows,decorations.markings,decorations.pathreplacing,quotes,fit,shapes,positioning}

\theoremstyle{remark}

\def\ca#1{\mathcal{#1}}

\title{Note on Min-$k$-Planar Drawings of Graphs}

\hideLIPIcs\nolinenumbers

\author{Petr Hlin\v{e}n\'y}{Masaryk University, Brno, Czech Republic}{hlineny@fi.muni.cz}{https://orcid.org/0000-0003-2125-1514}{}
\author{Lili K\"{o}dm\"{o}n}{E\"{o}tv\"{o}s Lor\'{a}nd University, Hungary}{kodmonlili@student.elte.hu}{https://orcid.org/0000-0002-6320-0152}{}
\authorrunning{P.~Hlin\v{e}n\'y and L.~K\"{o}dm\"{o}n}
\Copyright{Petr Hlin\v{e}n\'y and Lili K\"{o}dm\"{o}n}


\funding{Lili K\"{o}dm\"{o}n has been supported by the Ministry of Innovation and Technology of Hungary from the National Research,
 Development and Innovation Fund, financed under the ELTE TKP 2021-NKTA-62 funding scheme.}

\begin{document}
\maketitle              
\begin{abstract}
The $k$-planar graphs, which are (usually with small values of $k$ such as $1,2,3$) subject to recent intense research,
admit a drawing in which edges are allowed to cross, but each one edge is allowed to carry at most $k$ crossings.
In recently introduced~[Binucci et al., GD~2023] min-$k$-planar drawings of graphs, edges may possibly carry more than $k$ crossings,
but in any two crossing edges, at least one of the two must have at most $k$ crossings.
In both concepts, one may consider general drawings or a popular restricted concept of drawings called {\em simple}.
In a simple drawing, every two edges are allowed to cross at most once, and any two edges which share a vertex are forbidden to cross.

While, regarding the former concept, it is for $k\leq3$ known (but perhaps not widely known) that every general $k$-planar graph 
admits a \mbox{simple $k$-planar} drawing and this ceases to be true for any~$k\geq4$, 
the difference between general and simple drawings in the latter concept is more striking.
We prove that there exist graphs with a min-$2$-planar drawing, or with a min-$3$-planar drawing avoiding crossings of adjacent edges,
which have {\em no simple} min-$k$-planar~drawings for arbitrarily large fixed~$k$.

\end{abstract}

\section{Introduction}

One of the most popular current research directions in graph drawing is going ``beyond planarity''
\cite{DBLP:journals/csur/DidimoLM19,DBLP:books/sp/20/HT2020}.
This somewhat broad direction can be described as considering drawings of graphs in which
edges may cross, but the overall pattern of edge crossings is restricted, usually in a local setting.
Some of the earliest examples are $1$-planar graphs (every edge may have at most $1$ crossing),
and popular extensions nowadays include many families among which we, for example, mention
$k$-planar- and fan-planar graphs, or $k$-quasiplanar graphs.
These diverse classes often share some nice properties of planar graphs,
such as having few edges (e.g., for $1$-planar graphs~\cite{https://doi.org/10.1002/mana.19861250122},
for $k$-planar graphs~\cite{DBLP:journals/combinatorica/PachT97},
and for $3$-quasiplanar graphs~\cite{DBLP:journals/combinatorica/AgarwalAPPS97}).
However, they may differ greatly from planar graphs in other respects; for instance, recognizing $1$-planar graphs is NP-complete~\cite{DBLP:journals/algorithmica/GrigorievB07,DBLP:conf/gd/KorzhikM08}.

We refer to \Cref{sec:basics} for a precise definition of a drawing of a graph,
and of simple and min-$k$-planar drawings.
Briefly, by a {\em drawing} of a graph (here exclusively in the plane) we mean a topological representation
in which edges (as arcs) join their end vertices (as points) and avoid passing through other vertices.
Furthermore, every two distinct edges intersect finitely many times, and at most two edges intersect
in one point except when it is their common end vertex.
In a {\em simple drawing}, we additionally require that every two distinct edges intersect in at most one point -- a crossing or a common end.
(However, one has to be careful with this simplified definition when considering a simple drawing of a non-simple graph;
then two parallel edges share two common ends.)

It is very common that research papers assume only simple drawings, for convenience,
but it is sometimes not quite clear whether this assumption is made ``without loss of generality,'' or whether
it is a significant restriction on the kinds of drawings considered.
For instance, when studying the {\em crossing number} of a graph (the minimum total number of crossings over all drawings), 
one can quite easily restrict attention to simple drawings without loss of generality, but such restriction is not possible
when studying the odd crossing number~\cite{DBLP:journals/dcg/PelsmajerSS08}.
On a different note, so-called fan-planar graphs can be assumed to have a simple fan-planar drawing without loss of generality,
but the proof~\cite{DBLP:journals/jgaa/KlemzKRS23} is highly nontrivial.

\smallskip
Consider {\em$k$-planar graphs}, which are graphs admitting a drawing in which no edge carries more than $k$ crossings.
(The same concept is also known as the local crossing number \cite{schaefer2017crossing} or the crossing parameter \cite{DBLP:journals/algorithmica/GrigorievB07}.)
The seminal paper of Pach and T\'oth~\cite{DBLP:journals/combinatorica/PachT97} explicitly requires simple $k$-planar drawings,
while, e.g., ``algorithmic'' papers Grigoriev and Bodlaender~\cite{DBLP:journals/algorithmica/GrigorievB07} 
and Korzhik and Mohar~\cite{DBLP:conf/gd/KorzhikM08} deal with general $k$-planar drawings
(in fact, \cite{DBLP:conf/gd/KorzhikM08} mentions that any $1$-planar drawing can be turned into a simple one)
and Ackerman~\cite{DBLP:journals/comgeo/Ackerman19} distinguishes the cases.
A recent survey on beyond planarity \cite{DBLP:journals/csur/DidimoLM19} unfortunately does not explicitly address this issue,
and it mostly only implicitly restricts to results about simple drawings in this respect.

To illustrate the potential problem points (of unwary mixing general and simple $k$-planar drawings in research), we mention, e.g., 
\cite{DBLP:conf/gd/Bekos0R16} which formulates results about general $3$-planar drawings, but importantly uses a lemma of \cite{DBLP:journals/combinatorica/PachT97} which, in unmentioned fact, relies on the assumption of a simple drawing.
In this particular case of \cite{DBLP:conf/gd/Bekos0R16}, as well as in other papers which deal with $k$-planar drawings for only $k\leq3$, there is no reall problem since every general $k$-planar graph for $k\leq3$ admits a {simple $k$-planar} drawing, as shown already by Pach et al.~\cite{DBLP:journals/dcg/PachRTT06}.
On the other hand, for every $k\geq4$ there exist $k$-planar graphs which have no simple $k$-planar drawing, e.g., Schaefer~\cite[Chapter~7]{schaefer2017crossing}.

\smallskip

Concerning the new related concept of {\em min-$k$-planar graphs}, which are graphs admitting a drawing in which every pair of crossing
edges has one of the two edges with at most $k$ crossings; the introductory paper by Binucci et al.~\cite{DBLP:conf/gd/BinucciBBDDHKLMT23,DBLP:journals/corr/abs-2308-13401}
requires simple drawings by the definition.
However, min-$k$-planar drawings may also be understood in the general (non-simple) setting,
and the difference between the general and the simple settings is much more striking than in the case of $k$-planar graphs.

Namely, we prove (\Cref{thm:nosimple}) that for arbitrarily large fixed $k$ there exist graphs that are
min-$2$-planar without restricting to simple drawings, but which have {\em no simple} min-$k$-planar drawing.
Alternatively, counterexample graphs with a min-$3$-planar drawing in which no two edges sharing a common vertex cross can also be constructed.
In other words, the concepts of simple and general min-$k$-planar drawings \emph{always significantly differ}, except in the trivial case of $k=1$
(in which we can easily simplify any min-$1$-planar drawing, cf.~\Cref{prop:tosimple}),
and they differ for~$k\geq3$ even if we forbid general drawings in which two edges with a common vertex cross.

In the course of proving this result, we develop a technical tool (\Cref{lem:fullfram}) which suitably constrains
possible min-$k$-planar drawings of graphs within a rigid ``frame'', and we suggest this tool can be useful in further
research of the properties of min-$k$-planar graphs.

\section{Min-$k$-planar Drawings and Graphs}
\label{sec:basics}

We consider general finite undirected graphs (with possible parallel edges or loops), and say that a graph is {\em simple} if it has no parallel edges and no loops.

A \emph{drawing} $\mathcal{G}$ of a graph $G$ in the Euclidean plane $\mathbb{R}^2$ is a function that maps each vertex $v \in V(G)$ to a distinct point $\mathcal{G}(v) \in \mathbb R^2$ 
and each edge $e=uv \in E(G)$ to a simple curve (non-self-intersecting Jordan arc) $\mathcal{G}(e) \subset \mathbb R^2$ with the ends $\mathcal{G}(u)$ and $\mathcal{G}(v)$.
We require that $\ca G(e)$ is disjoint from $\ca G(w)$ for all~$w\in V(G)\setminus\{u,v\}$
and that $\ca G(e)\cap\ca G(f)$ is finite for all $e\not=f\in E(G)$.
In a slight abuse of notation we identify a vertex $v$ with its image $\mathcal{G}(v)$ and an edge $e$ with~$\mathcal{G}(e)$.
An intersection (possibly tangential) of two edges $e$ and $f$ other than a common end vertex,
in other words a point~$x\in\ca G(e)\cap\ca G(f)$, is called a \emph{crossing} (of $e$ and~$f$),
and the pair $e,f$ is said to {\em cross} in~$x$.
A {\em planar drawing} is a drawing with no crossings.

For further concepts, the definition of a `crossing' needs a careful clarification.
If $k>2$ edges together cross in one point, then this situation, strictly, counts as $k\choose2$ crossings of pairs of these edges
(see, e.g., \cite{DBLP:journals/siamcomp/CabelloM13}).
Since a slight local perturbation of the involved edges avoids common crossings of triples without adding any new crossing,
for our purposes, we may simply discard the possibility that $k>2$ edges together cross in one point, which is a common approach in literature.
Therefore, without loss of generality, we require that, in every drawing, {\em no three edges cross} in the same~point, unless stated otherwise.

One may, likewise, deal with possible tangential crossings in a drawing, which can be discarded by a slight local perturbation as well,
and so exclude tangential crossings by the definition (as many papers do for convenience), but here we stay on the more general side
and do not exclude them.

In a \emph{simple drawing} $\mathcal{G}$ (otherwise known as a {\em good drawing}~\cite{cd98262b-e086-3d8d-b726-cd66f05da520,schaefer2017crossing}, 
but we wish to adhere to the recent terminology), crossings are allowed, but (again)
no three edges cross in the same point, no two edges have more than one crossing in common, 
and no two adjacent edges (i.e., with a common end vertex) cross.
Hence, in our setting, a simple drawing $\mathcal{G}$ is defined such that, for every $e\not=f\in E(G)$,
we have $|\mathcal{G}(e)\cap\mathcal{G}(f)|\leq1$, except that $|\mathcal{G}(e)\cap\mathcal{G}(f)|=2$ when $e$ and $f$ are parallel edges
(the latter is irrelevant for simple~graphs).

A drawing $\mathcal{G}$ is {\em$k$-planar} if no edge contains more than $k$ crossings.
A drawing $\mathcal{G}$ is {\em min-$k$-planar} if, for every two crossing edges $e$ and $f$ in $\mathcal{G}$,
one (or both) of $e,f$ has no more than $k$ crossings.
If an edge $e$ has more than $k$ crossings in a min-$k$-planar drawing~$\mathcal{G}$, then~$e$ is called {\em heavy} in $\mathcal{G}$
(hence two heavy edges cannot cross each other in a min-$k$-planar~drawing).

A~graph $G$ is {\em min-$k$-planar} if $G$ admits a min-$k$-planar drawing.
Moreover, a graph $G$ is {\em simply min-$k$-planar} if $G$ admits a min-$k$-planar drawing $\mathcal{G}$ such that $\mathcal{G}$ is a simple drawing.
Since some papers shortly speak about \mbox{min-$k$-planar} graphs while requiring simple drawings, to avoid further confusion,
we will call min-$k$-planar graph without the additional requirement of a simple drawing as {\em general min-$k$-planar}.

A careful distinction between general min-$k$-planar graphs and simply min-$k$-planar graphs is indeed necessary
for all $k>1$, as we are going to prove here.

\begin{theorem}[Proof in \Cref{sec:details}]\label{thm:nosimple}
\begin{itemize}
\item[a)] For every $k\geq2$, there exists a simple graph $H_k$ which is general min-$2$-planar, but $H_k$ has {\em no simple} min-$k$-planar drawing.
\item[b)] Moreover, for all $k\geq3$, there exists a graph $H'_k$ which has a general min-$3$-planar drawing in which no two adjacent edges cross,
but, again, $H'_k$ has {\em no simple} min-$k$-planar drawing.
\end{itemize}
\end{theorem}

To complement \Cref{thm:nosimple}, we resolve the remaining trivial cases in \Cref{prop:tosimple}.

\begin{proposition}\label{prop:tosimple}
\begin{itemize}
\item[\it a)] Every general min-$1$-planar graph admits a simple min-$1$-planar drawing (hence is simply min-$1$-planar).
\item[\it b)] Every graph with a min-$2$-planar drawing in which no two adjacent edges cross also admits a simple min-$2$-planar drawing.
\end{itemize}
\end{proposition}

\begin{figure}[th]
  \centering
  \begin{minipage}[b]{0.4\textwidth}
\begin{tikzpicture}[xscale=1.2, yscale=0.9]\small
\tikzstyle{every path}=[draw, semithick]
\draw [blue!60!black] plot [smooth] coordinates {(-1.5,-0.1) (-0.7, 0.35) (0,0.5) (0.7, 0.35) (1.5,-0.1) (2.5,-1)};
\draw [red!60!black] plot [smooth] coordinates {(-1.5,-0.1) (-0.7, -0.55) (0,-0.7) (0.7, -0.55) (1.5,-0.1) (2.5,1)};
\draw [red!60!black, dashed] plot [smooth] coordinates {(-1.5,-0.1) (-0.7, -0.45) (0,-0.6) (0.7, -0.45) (1.5,-0.1)};
\tikzstyle{every node}=[draw, color=black, shape=circle, inner sep=1pt, fill=black]
\node[black] at (-1.5,-0.1) {};
\node[black, fill=none] at (1.5,-0.1) {};
\tikzstyle{every node}=[draw=none, fill=none]
\node[black] at (0,-0.95){$e_1$};
\node[black] at (0,0.75) {$f_1$};
\node[black] at (2.2,0.9) {$e$};
\node[black] at (2.2,-1.05) {$f$};
\node[black] at (-1.5,0.15) {$x$};
\node[black] at (1.5,0.2) {$y$};
\node[black] at (0,-0.3) {$e_1'$};
\draw[thin] (0,0.1) to[bend right=22] ++(1,1);
\draw[thin] (1.6,-0.7) to[bend right=14] ++(1,1);
\end{tikzpicture}
  \end{minipage}
  \hfill
  \begin{minipage}[b]{0.4\textwidth}
\begin{tikzpicture}[xscale=1.2, yscale=0.95]\small
\tikzstyle{every path}=[draw, semithick]
\draw [black, thin,dashed] plot [smooth] coordinates {(-1.5,-0.1) (-0.7, 0.35) (0,0.5) (0.7, 0.35) (1.5,-0.1) (2.5,-1)};
\draw [black, thin,dashed] plot [smooth] coordinates {(-1.5,-0.1) (-0.7, -0.55) (0,-0.7) (0.7, -0.55) (1.5,-0.1) (2.5,1)};
\draw [blue!60!black] plot [smooth] coordinates {(-1.5,-0.1) (-0.7, -0.65) (0,-0.8) (0.7, -0.65) (1.5,-0.35) (2.5,-1.1)};
\draw [red!60!black] plot [smooth] coordinates {(-1.5,-0.1) (-0.7, -0.45) (0,-0.6) (0.7, -0.45) (1.5,0) (2.5,1.1)};
\tikzstyle{every node}=[draw, color=black, shape=circle, inner sep=1pt, fill=black]
\node[black] at (-1.5,-0.1) {};
\node[black, fill=none] at (1.5,-0.1) {};
\tikzstyle{every node}=[draw=none, fill=none]
\node[black] at (-1.5,0.25) {$x$};
\node[black] at (1.5,0.25) {$y$};
\node[black] at (2.2,1) {$e''$};
\node[black] at (2.2,-1.2) {$f''$};
\draw[thin] (0,0.1) to[bend right=22] (1,1);
\draw[thin] (1.6,-0.7) to[bend right=14] ++(1,1);
\end{tikzpicture}
  \end{minipage}
\caption{An illustration of \Cref{prop:tosimple}\,a); up to symmetry between $e$ and $f$ (with common end vertex $x$), the edge $e$ carries no other crossing than the point $y$, and so one can draw an uncrossed arc $e'_1$ tightly along the segment $e_1\subseteq e$ from $x$ to $y$.
When redrawing from $e$ to $e''$ and from $f$ to $f''$ (using $e'_1$), the crossing at $y$ is eliminated and no new crossings are added to any edge in the picture.}
\label{fig:simplif1}
\end{figure}
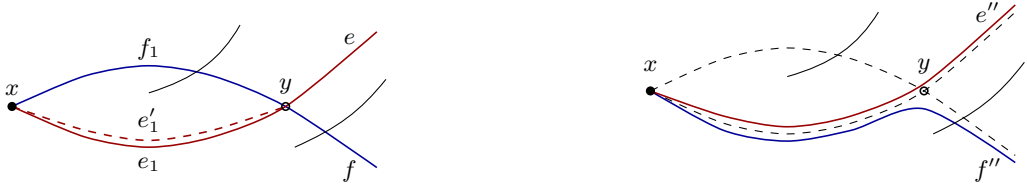

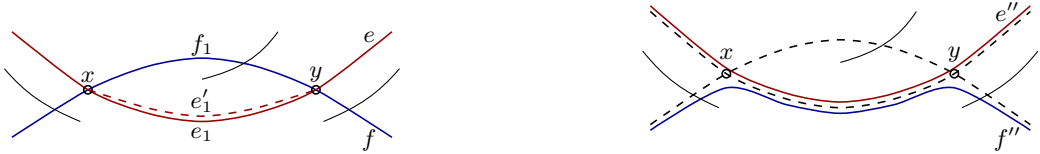
\begin{figure}[!tbp]
  \centering
  \begin{minipage}[b]{0.4\textwidth}
\begin{tikzpicture}[xscale=1, yscale=0.7]\small
\tikzstyle{every path}=[draw, semithick]
\draw [blue!60!black] plot [smooth] coordinates {(-2.5,-1) (-1.5,-0.1) (-0.7, 0.35) (0,0.5) (0.7, 0.35) (1.5,-0.1) (2.5,-1)};
\draw [red!60!black] plot [smooth] coordinates {(-2.5,1) (-1.5,-0.1) (-0.7, -0.55) (0,-0.7) (0.7, -0.55) (1.5,-0.1) (2.5,1)};
\draw [red!60!black, dashed] plot [smooth] coordinates {(-1.5,-0.1) (-0.7, -0.45) (0,-0.6) (0.7, -0.45) (1.5,-0.1)};
\tikzstyle{every node}=[draw, color=black, shape=circle, inner sep=1.1pt, fill=none]
\node[black] at (-1.5,-0.1) {};
\node[black] at (1.5,-0.1) {};
\tikzstyle{every node}=[draw=none, fill=none]
\node[black] at (0,-0.95){$e_1$};
\node[black] at (0,0.75) {$f_1$};
\node[black] at (2.2,0.9) {$e$};
\node[black] at (2.2,-1.05) {$f$};
\node[black] at (-1.5,0.15) {$x$};
\node[black] at (1.5,0.2) {$y$};
\node[black] at (0,-0.3) {$e_1'$};
\draw[thin] (0,0.1) to[bend right=22] (1,1);
\draw[thin] (1.6,-0.7) to[bend right=14] ++(1,1);
\draw[thin] (-1.6,-0.7) to[bend left=14] ++(-1,1);
\end{tikzpicture}
  \end{minipage}
  \hfill
  \begin{minipage}[b]{0.4\textwidth}
\begin{tikzpicture}[xscale=1, yscale=0.75]\small
\tikzstyle{every path}=[draw, semithick]
\draw [black, dashed] plot [smooth] coordinates {(-2.5,-1) (-1.5,-0.1) (-0.7, 0.35) (0,0.5) (0.7, 0.35) (1.5,-0.1) (2.5,-1)};
\draw [black, dashed] plot [smooth] coordinates {(-2.5,1) (-1.5,-0.1) (-0.7, -0.55) (0,-0.7) (0.7, -0.55) (1.5,-0.1) (2.5,1)};
\draw [blue!60!black] plot [smooth] coordinates {(-2.5,-1.1) (-1.5,-0.35) (-0.7, -0.65) (0,-0.8) (0.7, -0.65) (1.5,-0.35) (2.5,-1.1)};
\draw [red!60!black] plot [smooth] coordinates {(-2.5,1.1) (-1.5,0) (-0.7, -0.45) (0,-0.6) (0.7, -0.45) (1.5,0) (2.5,1.1)};
\tikzstyle{every node}=[draw, color=black, shape=circle, inner sep=1.1pt, fill=none]
\node[black] at (-1.5,-0.1) {};
\node[black] at (1.5,-0.1) {};
\tikzstyle{every node}=[draw=none, fill=none]
\node[black] at (-1.5,0.25) {$x$};
\node[black] at (1.5,0.25) {$y$};
\node[black] at (2.2,1) {$e''$};
\node[black] at (2.2,-1.25) {$f''$};
\draw[thin] (0,0.1) to[bend right=22] (1,1);
\draw[thin] (1.6,-0.7) to[bend right=14] ++(1,1);
\draw[thin] (-1.6,-0.7) to[bend left=14] ++(-1,1);
\end{tikzpicture}
  \end{minipage}
\caption{An illustration of \Cref{prop:tosimple}\,b); now we have two crossings $x$ and $y$ of the same pair $e$ and $f$ of edges,
and there are no more crossings on $e$.
Similarly to \Cref{fig:simplif1}, when redrawing from $e$ to $e''$ and from $f$ to $f''$ (using the uncrossed arc $e'_1$), at least one of the crossings at $x$ or $y$ is eliminated and no new crossings are added to any edge in the picture.}
\label{fig:simplif2}
\end{figure}

\begin{proof}
a) Let $G$ be a general min-$1$-planar graph and take a min-$1$-planar drawing $\mathcal{G}$ of~$G$.
For simplicity, in this proof, we write $e$ (an edge) also for the point set $\mathcal{G}(e)$ in the drawing.
We may assume that $\mathcal{G}$ minimizes the number of edge pairs which violate the simplicity of the drawing,
i.e., edge pairs which share a vertex and cross, since a pair cannot cross twice in a min-$1$-planar drawing.
Let $e,f\in E(G)$ be such a violating edge pair in $\mathcal{G}$;
hence the intersection of $e$ and $f$ contains two distinct points $x,y\in e\cap f$ where $x$ is a vertex and $y$ a crossing.
Moreover, we may assume up to symmetry that $e$ has no crossing other than with~$f$.

Let us denote by $e_1$ and $f_1$ the subarcs of $e$ and $f$, respectively, in $\mathcal{G}$ with the ends~$x,y$.
Since $e_1$ is internally crossing-free, there exists an arc $e_1'$ from $x$ to $y$ drawn sufficiently close to $e_1$
such that $e_1'$ is disjoint from $\mathcal{G}$ except in~$x,y$.
We replace $f$ in $\mathcal{G}$ with $f':=(f\setminus f_1)\cup e_1'$;
this new drawing $\mathcal{G}'$ of $G$ is again min-$1$-planar since no crossings of~$\mathcal{G}$ have been affected.
If $y$ is a tangent of the edges $e$ and $f'$ (not crossing transversely at~$y$), then a local perturbation
of $e$ around $y$ simply removes this tangential crossing, and so we eliminate one violating edge pair from~$\mathcal{G}$ in~$\mathcal{G}'$.
Otherwise, if $e$ and $f'$ cross transversely at~$y$, we instead replace $e$ with $e'':=(e\setminus e_1)\cup e_1'$ and $f$ with~$f'':=(f\setminus f_1)\cup e_1$ in~$\mathcal{G}$.
See \Cref{fig:simplif1}.
The new drawing $\mathcal{G}''$ of $G$ is min-$1$-planar, too, and $y$ is now a tangent of the edges $e''$ and $f''$, which can be eliminated as previously.

b) We let $\mathcal{G}$ be a min-$2$-planar drawing of a graph $G$ in which no two adjacent edges cross.
Again, we may assume that $\mathcal{G}$ minimizes the number of edge pairs which violate the simplicity of the drawing
-- these are now the pairs which mutually cross exactly twice.
Let $e,f\in E(G)$ be such a violating edge pair in $\mathcal{G}$;
hence the intersection of $e$ and $f$ contains two distinct points $x,y\in e\cap f$ which are crossings.

Since $\mathcal{G}$ is min-$2$-planar, up to symmetry, the edge $e$ has no other crossings than~$x,y$.
Denote by $e_1$ and $f_1$ the subarcs of $e$ and $f$, respectively, in $\mathcal{G}$ with the ends~$x,y$.
Since $e_1$ is internally crossing-free, there exists an arc $e_1'$ from $x$ to $y$ drawn sufficiently close to $e_1$
such that $e_1'$ is disjoint from $\mathcal{G}$ except in~$x,y$.
We define new arcs of the edges $e$ and $f$ in $\mathcal{G}$ as follows;
$e'=e$ and $f':=(f\setminus f_1)\cup e_1'$, and $e'':=(e\setminus e_1)\cup e_1'$ and $f'':=(f\setminus f_1)\cup e_1$.
Clearly, for at least one of the pairs $e',f'$ or $e'',f''$, some of the points $x,y$ is now a tangential crossing, which can be eliminated as previously.
See \Cref{fig:simplif2}.
The new drawing $\mathcal{G}'$ of $G$ is min-$2$-planar, too, and no new crossings have been created between any edge pairs.

Altogether, we have in each case decreased the number of violating edge pairs, contradicting our minimal choice of the drawing~$\mathcal{G}$.
\end{proof}

To prove \Cref{thm:nosimple}, we use the following intermediate result formulated for so-called anchored graphs,
which captures the essence of \Cref{thm:nosimple}.
An {\em anchored graph} is a pair $(G,A)$ where $A\subseteq V(G)$ is an ordered tuple of vertices.
An {\em anchored drawing} of $(G,A)$ in the unit disk $D\subseteq\mathbb{R}^2$ is a drawing $\mathcal{G}\subseteq D$ of $G$ 
such that $\mathcal{G}$ intersects the boundary of $D$ precisely in the points of $A$ (the {\em anchors}) in this clockwise order.
We naturally extend the adjective {\em anchored} to min-$k$-planar drawings.
We prove:

\begin{lemma}\label{lem:nosimplea}
For every $k\geq2$, there exists a simple anchored graph $(G_k,A_k)$ which has an anchored general min-$2$-planar drawing,
but $(G_k,A_k)$ has {\em no} anchored {\em simple} min-$k$-planar drawing.
Furthermore, for any $k\geq3$, there exists a simple anchored graph $(G'_k,A'_k)$ which has an anchored general min-$3$-planar drawing
in which no pair of adjacent edges cross, but $(G'_k,A'_k)$ has {\em no} anchored {\em simple} min-$k$-planar drawing.
\end{lemma}

\begin{proof}
We define the anchored graph $(G_k,A_k)$ as depicted in \Cref{fig:anch2no}.
$G_k$ is a disjoint union of two induced matchings $M_1$ and $M_2$ of $k+1$ edges each, and of two induced stars $S_1$~and~$S_2$,
where $S_1$ has center $a_2$ and $k+1$ leaves including $a_1,b_1$ and $S_2$ has center $c_2$ and $k+2$ leaves in the set $C_3\cup\{c_1\}$ where $|C_3|=k+1$.
The anchor set is $A_k=V(G_k)\setminus\{c_2\}$ ordered as in \Cref{fig:anch2no}.
\Cref{fig:anch2no} also shows an anchored min-$2$-planar drawing of $(G_k,A_k)$ (however, $a_1a_2$ crosses~$b_1a_2$ there).

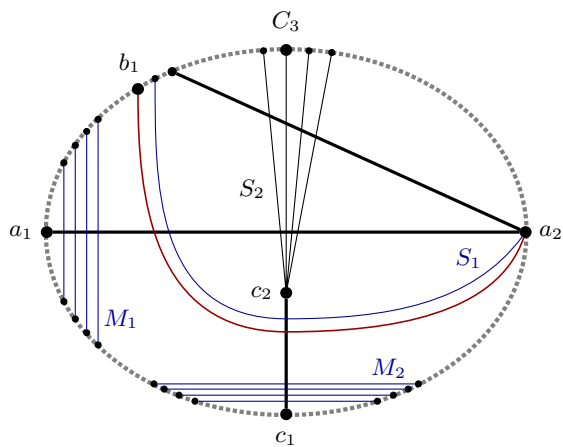
\begin{figure}[tb]
$$
\begin{tikzpicture}[xscale=1.5, yscale=1.15]\small
\draw[gray, ultra thick, densely dotted] (0,0) ellipse (60pt and 60pt);
\tikzstyle{every node}=[draw, color=black, thick, shape=circle, inner sep=1.3pt, fill=black]
\node[label=left:$a_1$] (a1) at (-2.1,0) {}; \node[label=right:$a_2$] (a2) at (2.1,0) {};
\node[label=below:$c_1$] (a3) at (0,-2.1) {}; \node[label=left:$c_2$] (c) at (0,-0.7) {};
\node[label=above:$\!\!b_1~$] (b1) at (-1.3,1.65) {};
\node[label=above:$C_3$] (w) at (0,2.1) {};
\tikzstyle{every node}=[draw, color=black, shape=circle, inner sep=0.8pt, fill=black]
\node (v1) at (-0.8,-1.95) {}; \node (v2) at (0.8,-1.95) {};
\node (v7) at (-0.94,-1.88) {}; \node (v8) at (0.94,-1.88) {};
\node (v3) at (-1.07,-1.81) {}; \node (v4) at (1.07,-1.81) {};
\node (v5) at (-1.16,-1.75) {}; \node (v6) at (1.16,-1.75) {};
\node[very thick] (t3) at (-1,1.85) {}; 
\node (u1) at (-1.95,-0.8) {}; \node (u2) at (-1.95,0.8) {};
\node (u3) at (-1.85,-1) {}; \node (u4) at (-1.85,1) {};
\node (u5) at (-1.75,-1.16) {}; \node (u6) at (-1.75,1.16) {};
\node (u7) at (-1.65,-1.3) {}; \node (u8) at (-1.65,1.3) {};
\draw[very thick] (a1)--(a2) (a3)--(c) (t3)--(a2);
\draw[blue!60!black] (v1)--(v2) (v3)--(v4) (v5)--(v6) (v7)--(v8) ;
\draw[blue!60!black] (u1)--(u2) (u3)--(u4) (u5)--(u6) (u7)--(u8) ;
\draw (c)--(w) (c)--(0.2,2.09) node{} (c)--(-0.2,2.09) node{} (c)--(0.4,2.07) node{} ;
\draw[blue!60!black] (-1.15,1.77) node{} to[out=270,in=180] (0,-1.0) to[out=0,in=240] (a2);
\draw[red!60!black,semithick] (b1) to[out=270,in=180] (0,-1.15) to[out=0,in=255] (a2);
\tikzstyle{every node}=[draw=none, fill=none]
\node[blue!60!black] at (-1.45,-1) {$M_1$}; \node[blue!60!black] at (0.9,-1.57) {$M_2$};
\node at (-0.3,0.5) {$S_2$}; \node[blue!60!black] at (1.6,-0.3) {$S_1$};
\end{tikzpicture}
\vspace*{-2ex}$$
\caption{An anchored general min-$2$-planar drawing of the anchored graph $(G_k,A_k)$ (here for $k=3$, with the heavy edges drawn thick) as used in the proof of \Cref{lem:nosimplea}.
	With general $k\geq2$, at least $k-1$ of the edges of the star $S_1$ with center $a_2$ follow the depicted route of the red
	edge~$a_2b_1$, while crossing the edge~$a_2a_1$.}
\label{fig:anch2no}
\end{figure}

Assume, for a contradiction, that there exists an anchored simple min-$k$-planar drawing $\mathcal{G}$ of $(G_k,A_k)$.
By Jordan Curve Theorem, the edge $a_1a_2$ has to cross all $k+1$ edges of $M_1$, and so $a_1a_2$ is heavy in $\mathcal{G}$.
If the edge $c_1c_2$ was crossing $a_1a_2$, then, again using Jordan Curve Theorem,
some edge of $M_2$ would cross twice with $a_1a_2$ -- which is not simple,
or $c_1c_2$ would cross all edges of $M_2$ and be heavy as well -- which contradicts $\mathcal{G}$ being min-$2$-planar.
Therefore, all $k+1$ edges from $c_2$ to $C_3$ have to cross $a_1a_2$ and are non-heavy.

Consequently, none of the edges between $c_2$ and $C_3$ can cross all $k+1$ edges of $S_1$.
By Jordan Curve Theorem, hence, some edge of $S_1$, say $b_1a_2$ (as in the picture) crosses $c_1c_2$.
However, again by Jordan Curve Theorem, from this crossing point to $b_1$, the edge $b_1a_2$ has to cross the edge $a_1a_2$,
contradicting the assumption that $\mathcal{G}$ is simple.

\begin{figure}[tb]
$$
\begin{tikzpicture}[xscale=1.5, yscale=1.15]\small
\draw[gray, ultra thick, densely dotted] (0,0) ellipse (60pt and 60pt);
\tikzstyle{every node}=[draw, color=black, thick, shape=circle, inner sep=1.3pt, fill=black]
\node[label=left:$a_1$] (a1) at (-2.1,0) {}; \node[label=right:$a_2$] (a2) at (2.1,0) {};
\node[label=below:$c_1$] (a3) at (0,-2.1) {}; \node[label=left:$c_2$] (c) at (0,-0.85) {};
\node[label=above:$b_1~$] (b1) at (-1.2,1.72) {}; \node[label=above:$~b_2$] (b2) at (1.2,1.72) {};
\node[label=above:$c_3$] (w) at (0,2.1) {};
\tikzstyle{every node}=[draw, color=black, shape=circle, inner sep=0.8pt, fill=black]
\node (v1) at (-0.8,-1.95) {}; \node (v2) at (0.8,-1.95) {};
\node (v7) at (-0.91,-1.9) {}; \node (v8) at (0.91,-1.9) {};
\node (v3) at (-1,-1.85) {}; \node (v4) at (1,-1.85) {};
\node (v9) at (-1.08,-1.8) {}; \node (v10) at (1.08,-1.8) {};
\node (v5) at (-1.16,-1.75) {}; \node (v6) at (1.16,-1.75) {};
\node (t7) at (-0.8,1.95) {}; \node (t8) at (0.8,1.95) {};
\node (t3) at (-1,1.85) {}; \node (t4) at (1,1.85) {};
\node (t9) at (-1.08,1.8) {}; \node (t10) at (1.08,1.8) {};
\node (u1) at (-1.95,-0.8) {}; \node (u2) at (-1.95,0.8) {};
\node (u3) at (-1.85,-1) {}; \node (u4) at (-1.85,1) {};
\node (u5) at (-1.75,-1.16) {}; \node (u6) at (-1.75,1.16) {};
\node (u7) at (-1.65,-1.3) {}; \node (u8) at (-1.65,1.3) {};
\node (u9) at (-1.55,-1.42) {}; \node (u10) at (-1.55,1.42) {};
\draw[very thick] (a1)--(a2) (a3)--(c);
\draw[blue!60!black] (v1)--(v2) (v3)--(v4) (v5)--(v6) (v7)--(v8) (v9)--(v10) ;
\draw[blue!60!black] (u1)--(u2) (u3)--(u4) (u5)--(u6) (u7)--(u8) (u9)--(u10) ;
\draw[blue!60!black]  (t3)--(t4)  (t7)--(t8) ;
\draw[blue!60!black] (t9) to[out=270,in=180] (0,-1.22) to[out=0,in=270] (t10);
\draw (c)--(w) ;
\draw[red!60!black,semithick] (b1) to[out=270,in=180] (0,-1.4) to[out=0,in=270] (b2);
\tikzstyle{every node}=[draw=none, fill=none]
\node[blue!60!black] at (-1.35,-1) {$M_1$};
\node[blue!60!black] at (0.9,-1.55) {$M_2$}; \node[blue!60!black] at (0.9,1.6) {$M_3$};
\end{tikzpicture}
\vspace*{-2ex}$$
\caption{An anchored general min-$3$-planar drawing of the anchored graph $(G'_k,A'_k)$ (here for $k=4$, with the heavy edges drawn thick) as used in the proof of \Cref{lem:nosimplea}.
	For general $k\geq3$, at least~$k-2$ of the edges of the matching $M_3$ follow the route of the red edge $b_1b_2$.
	These edges cross the heavy edge $a_1a_2$ twice, and so the drawing is neither simple nor general min-$2$-planar.}
\label{fig:anch3no}
\end{figure}
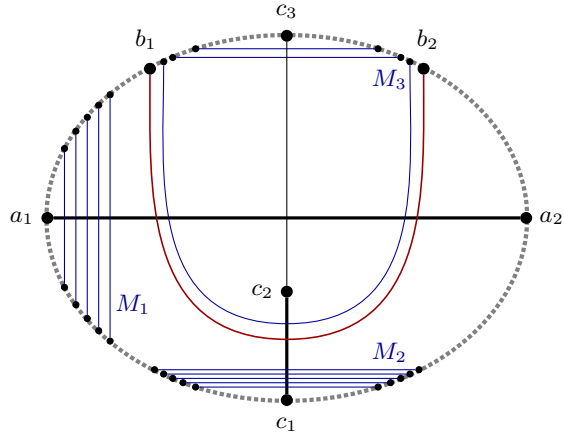

\smallskip
For $k\geq3$, we define the anchored graph $(G'_k,A'_k)$ as depicted (with a min-$3$-planar drawing) in \Cref{fig:anch3no}.
The definition is analogous to that of $(G_k,A_k)$, except that the star $S_1$ is now replaced with an induced matching
$M_3$ of $k+1$ edges (including the edge $a_1a_2$), and the star $S_2$ is replaced with a path $(c_1,c_2,c_3)$ of length two.

A proof that $(G'_k,A'_k)$ has no anchored simple min-$k$-planar drawing starts with the same steps as the previous one.
We get that the edge $a_1a_2$ is heavy and not crossing $c_1c_2$, and so the edge $c_2c_3$ has to cross $a_1a_2$ and is non-heavy.
Then each of the $k$ edges of $M_3$ has to cross $c_2c_3$ or $c_1c_2$, but not all may cross $c_2c_3$ which would become heavy.
Consequently, some edge, say $b_1b_2\in E(M_3)$, crosses the edge $c_1c_2$.
Since $c_1c_2$ is separated from $V(M_3)$ by $a_1a_2$, by Jordan Curve Theorem, the edge $b_1b_2$ has to cross $a_1a_2$ twice,
again a contradiction.
\end{proof}

\section{Technical Proofs}
\label{sec:details}

In this section we give the technical details leading to the proof of \Cref{thm:nosimple}.
In a nutshell, we are going to construct a ``frame graph'' which enforces a predefined anchored subdrawing 
of a given anchored graph, and then we apply this construction to the graphs of \Cref{lem:nosimplea}.

A {\em$t$-amplification} of a graph $G$ is the graph obtained from $G$ by replacing every edge
$e=xy\in E(G)$ with a new collection of $t$ pairwise internally disjoint paths of length~$2$ from $x$ to~$y$.
Observe that any $t$-amplification of a planar graph $G$ is again planar.

\begin{lemma}\label{lem:amplify}
For every planar graph $G$, there exists an integer $t=t(G,k,w)$, depending on $G$ and integers $k$ and~$w$, such that the following holds.
Every general min-$k$-planar drawing $\mathcal{G}$ of the $t$-amplification $G^t$ of $G$ contains a planar subdrawing $\mathcal{G}'$
which is isomorphic to the $w$-amplification $G^w$ of~$G$.
\end{lemma}

\begin{proof}
Let the length-$2$ paths (with ends in $V(G)$\,) from the definition of the $t$-amplification $G^t$ be called {\em double edges} of $G^t$.
We start with an easy claim:
\begin{itemize}
\item[(*)] Let $a_1,a_2,b_1,b_2\in V(G)$ (possibly $a_1=a_2$ or $b_1=b_2$), and let $D_1$ and $D_2$ be sets of double edges
of $G^t$ where those of $D_i$ have ends $a_i$ and $b_i$, such that each double edge of $D_1$ crosses each one of $D_2$ in the drawing $\mathcal{G}$.
If $|D_1|,|D_2|\geq2k+1$, then $\mathcal{G}$ is not min-$k$-planar.
\end{itemize}
Indeed, since every double edge $P\in D_1$ crosses all of $D_2$, one of the two edges of $P$ has at least $k+1$ crossings and is heavy.
The same holds for every $P'\in D_2$.
Without loss of generality, assume now $|D_1|=|D_2|=2k+1$.
If every crossing occurring between $D_1$ and $D_2$ involved a non-heavy edge, then the total number of these crossings
would be bounded (counted along all non-heavy edges) by at most $|D_1|\cdot k+|D_2|\cdot k=k(4k+2)<(2k+1)^2=|D_1|\cdot|D_2|$,
which is impossible.
Hence some two heavy edges cross and $\mathcal{G}$ is not min-$k$-planar.

\smallskip
We continue with a Ramsey-type argument.%
\footnote{One can get better estimates of $t$ using homotopy-based arguments, but that would not make our result stronger
and we prefer simplicity of brute-force Ramsey here.}
For any $f\in E(G)$ and any $t_1$, by Ramsey Theorem, there is a sufficiently large $t$ such that the following holds.
Among the $t$ double edges replacing $f$ in the $t$-amplification $G^t$, there exist $2(2k+1)$ pairwise crossing ones or
$t_1$ pairwise non-crossing ones.
If the former (crossing) case happens, then by (*) we get a contradiction to the assumption of this lemma, 
that $\mathcal{G}$ is min-$k$-planar.
Therefore, the non-crossing case happens, and we apply the same argument concurrently to all edges of~$G$.
This way we get a subdrawing $\mathcal{G}_1\subseteq\mathcal{G}$ of a $t_1$-amplification~$G^{t_1}$ of $G$ 
such that the collection of paths replacing any edge of $G$ is alone crossing-free in $\mathcal{G}_1$.

Likewise, for any pair of edges $f,f'\in E(G)$ and any $t_2$, by the bipartite Ramsey Theorem, 
there is a sufficiently large $t_1$ such that we get $t_2$ double edges replacing $f$ and another $t_2$ replacing $f'$ in~$G^{t_1}$,
which are pairwise noncrossing in $\mathcal{G}_1$, or we again get a contradiction via (*).
Applying this argument concurrently to disjoint pairs of edges of $G$, 
we obtain a subdrawing $\mathcal{G}_2\subseteq\mathcal{G}_1$ of a $t_2$-amplification of $G$.
We iterate this argument until we exhaust all pairs of edges of~$G$.
Starting from a sufficiently large $t$, the resulting drawing $\mathcal{G}'$ of this iterative process is
a $w$-amplification of $G$, and no two double edges cross in $\mathcal{G}'$, as desired.
\end{proof}

\begin{lemma}\label{lem:fullfram}
For any integers $a,k$ and simple graph $G$ with an ordered subset $A\subseteq V(G)$, \mbox{$|A|=a$}, 
there exists a simple anchored graph $(H,A)$ disjoint from $G$ except in the anchors~$A$, such that the following hold:
\begin{itemize}
\item[a)] $(H,A)$ has an anchored simple min-$1$-planar drawing.
\item[b)] In every general min-$k$-planar drawing $\mathcal{H}$ of $H\cup G$, the subdrawing $\mathcal{G}\subseteq\mathcal{H}$ of $G$
is (spherically) homeomorphic to an anchored drawing of $(G,A)$ or its mirror image.
\end{itemize}
\end{lemma}

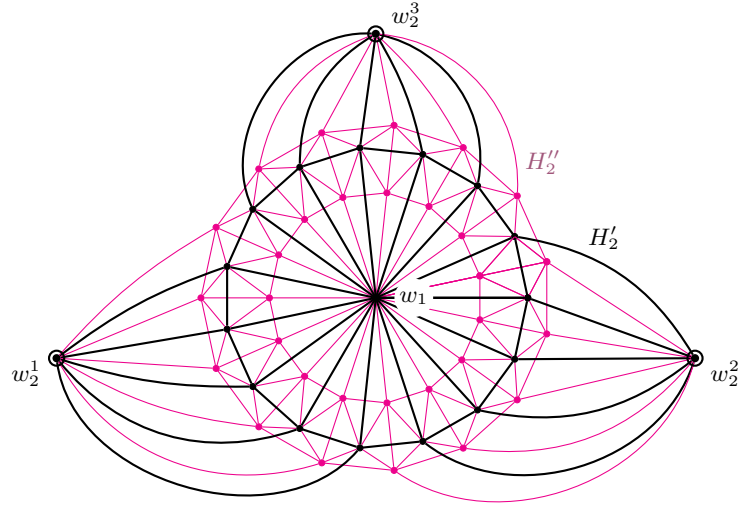
\begin{figure}[t]
$$
\begin{tikzpicture}[scale=1]\small
\coordinate (uu) at (0,0); \coordinate (uv) at (0,0);
\tikzstyle{every node}=[draw, color=black, shape=circle, inner sep=0.8pt, fill=black]
\foreach \u in {0,24,...,360} { \node (W\u) at (\u:2) {}; }
\coordinate (W384) at (W24);
\tikzstyle{every node}=[draw, color=magenta, shape=circle, inner sep=0.8pt, fill=magenta]
\tikzstyle{every path}=[draw, color=magenta]
\coordinate (uu) at (0,0);
\foreach \u in {12,36,...,372} {
	\node (M\u) at (\u:1.4) {};  \draw (0,0)--(M\u);
	\node (N\u) at (\u:2.3) {};  \draw (M\u)--(N\u);
	\draw (uu)--(M\u);  \coordinate (uu) at (M\u);
	\draw (uv)--(N\u);  \coordinate (uv) at (N\u);
	\pgfmathtruncatemacro{\um}{\u - 12}
	\pgfmathtruncatemacro{\up}{\u + 12}
	\draw (W\um)--(M\u)--(W\up) (W\um)--(N\u)--(W\up);
}
\tikzstyle{every node}=[draw, color=black, shape=circle, inner sep=0.8pt, fill=black]
\tikzstyle{every path}=[draw, color=black, thick]
\coordinate (uu) at (0,0);
\foreach \u in {0,24,...,360} {
	\draw (0,0)--(W\u);  \draw (uu)--(W\u);  \coordinate (uu) at (W\u);
}
\node (A) at (0,3.5) {};  \node (B) at (4.2,-0.8) {};  \node (C) at (-4.2,-0.8) {};
\draw (A)--(W96) (A) to[bend left=9] (W72) (A) to[bend left=44] (W48) (A) to[bend right=29] (W120) (A) to[bend right=59] (W144);
\draw (B)--(W0) (B)--(W336) (B) to[bend left=25] (W312) (B) to[bend left=55] (W288) (B) to[bend right=25] (W24);
\draw (C)--(W192) (C) to[bend right=9](W216) (C) to[bend right=35] (W240) (C) to[bend left=9] (W168) (C) to[bend right=66] (W264);
\tikzstyle{every path}=[draw, color=magenta]
\draw (A)--(N84) (A) to[bend left=9] (N60) (A) to[bend left=44] (N36) (A)--(N108) (A) to[bend right=29] (N132);
\draw (B)--(N12) (B) to[bend left=55] (N276) (B)--(N348) (B)--(N324) (B) to[bend left=29] (N300);
\draw (C)--(N180) (C) to[bend left=9] (N156) (C) to[bend right=39] (N252) (C) to[bend right=12] (N228) (C)--(N204);
\tikzstyle{every node}=[draw, color=black, shape=circle, inner sep=2pt, thick, fill=none]
\node at (A) {}; \node at (B) {}; \node at (C) {};
\tikzstyle{every node}=[draw=none, fill=none, shape=circle, inner sep=1pt]
\node[black] at (3,0.8) {$H_2'$};
\node[black,fill=white] at (0.5,0) {$w_1$};
\node[black] at (-4.6,-1) {$w_2^1$};
\node[black] at (4.6,-1) {$w_2^2$};
\node[black] at (0.4,3.75) {$w_2^3$};
\node[magenta!60!black] at (2.2,1.8) {$H_2''$};
\end{tikzpicture}
\vspace*{-2ex}$$
\caption{An illustration of the graph $H_2$ from the proof of \Cref{lem:fullfram}.
	The three black circled vertices are the designated anchors~$A$ (here~$A=\{w_2^1,w_2^2,w_2^3\}$).
	The subgraph $H_2'$ is in thick black and $H_2''$ in magenta colour.
	All magenta edges get $t$-amplified in the construction of~$H$, and consequently,
	the black edges forced to cross them will be heavy in any min-$k$-planar drawing of $H$ by \Cref{lem:amplify}.}
\label{fig:dualH2}
\end{figure}

\begin{proof}
For a start, we ignore all components of $G-A$ which attach to $A$ in at most one vertex;
their possible subdrawings can always be added homeomorphically and without further crossings to the rest of an anchored drawing of $(G,A)$.
Let $\ell$ be the maximum finite(!) distance in $G$ between a vertex of $A$ and a (reachable) vertex of $V(G)\setminus A$.
Let $d=a(2k\ell+2k+1)$.

We begin with the graph $H_0$ which is a double wheel of $d$ spokes, i.e., a graph made from the cycle $C_d$ by adding
two central vertices $w_1,w_2$ adjacent to all cycle vertices.
Let $H_0^*$ be the planar dual of $H_0$, and $H_1$ be constructed from $H_0\cup H_0^*$ by adding extra edges
$uv\in E(H_1)$ for every pair $u\in V(H_0)$ and $v\in V(H_0^*)$ such that $u$ is a vertex incident
to the face of (the unique planar drawing of) $H_0$ represented by $v$.
We construct the graph $H_2$ by splitting the central vertex $w_2$ into $a$ vertices, each incident with $2k\ell+2k+1$ consecutive
spokes of $H_0$ and with naturally corresponding extra edges of~$H_1$. See \Cref{fig:dualH2}.
The vertices split from $w_2$, in their clockwise order, define the anchor set~$A$.

Furthermore, let $H_2'\subseteq H_2$ be the subgraph formed by the original edges of $H_0$, and $H_2''=H_2\setminus E(H_2')$.
Observe one important property of planar $H_2''$; for every planar drawing of $H_2''$ and each edge $uv\in E(H_2')$,
the points $u$ and $v$ are separated by a cycle (e.g., the dual cycle of $u\not=w_2$) in~$H_2''$.
Let now $t=t(H_2'',k,k+1)$ be as in \Cref{lem:amplify}.
We construct $H:=H_2'\cup H_2^t$ where $H_2^t$ is the $t$-amplification of~$H_2''$.
Obviously, $(H,A)$ has an anchored simple min-$1$-planar drawing, e.g., one following \Cref{fig:dualH2}.

\smallskip
Consider now a min-$k$-planar drawing $\mathcal{H}$ of $H\cup G$, 
and denote by $\mathcal{H}_0\subseteq\mathcal{H}$ the subdrawing of~$H_2^t$.
By \Cref{lem:amplify}, $\mathcal{H}_0$ contains a planar subdrawing of a $(k+1)$-amplification of $H_2''$.
By the mentioned property of $H_2''$, every edge of $H_2'$ thus has to cross at least $k+1$ edges of $H_2^t$ and is heavy.
In particular, the subdrawing of $H_2'$ within $\mathcal{H}$ is planar (and so has to look as in \Cref{fig:dualH2}).

We pick an anchor vertex $a\in A$ and observe that $a$ has $2k\ell+2k+1$ disjoint length-$2$ paths to~$w_1$ in $H_2'$.
Let $P_a$ denote the ``middle'' one of them.
Assume that an edge $f=xy\in E(G)$ crosses $P_a$ in the drawing $\mathcal{H}$.
Then each end $x$ or $y$ has distance at most $\ell$ in $G$ to a vertex $b\in A\setminus\{a\}$ by our choice of $\ell$,
and so $\leq\ell+1$ edges of a path from $f$ (and including $f$) to~$b$ have to cross together $k\ell+k+1$ paths of $H_2'$ 
between $a$ and~$w_1$, by Jordan Curve Theorem.
Therefore, some edge of this path of $G$ must cross at least $k+1$ edges of $H_2'$, and we get pair(s) of
crossing heavy edges, which is a contradiction.

Therefore, in $\mathcal{H}$ there are $H_2'$-paths from each anchor vertex $a\in A$ to central $w_1$ which are not crossed by
any edge of the subdrawing of~$G$.
This means that the subdrawing of $G$ within $\mathcal{H}$ can be homeomorphically deformed in the sphere so that the vertices
of $A$ will be drawn on a disk boundary and the rest inside, as required by an anchored drawing.
The correct cyclic order (up to a mirror image) of the anchors $A$ on the disk boundary is ensured by the cycle
of $H_2'$ on the neighbours of~$w_1$ (the rim cycle of the starting wheel).
\end{proof}

We are now ready to finish the proof of the main result.

\begin{proof}[Proof of \Cref{thm:nosimple}]
We take the anchored graph $(G_k,A_k)$, $k\geq2$, and the anchored graph $(G'_k,A'_k)$, $k \geq 3$, from \Cref{lem:nosimplea} for (a) and (b), respectively,
and plug them into \Cref{lem:fullfram} as $G$.
The rest is an immediate consequence of the previous statements.
\end{proof}

\section{Conclusions}
\label{sec:conclu}

Many papers in the graph drawing area deal with only simple drawings, either as a convenient simplification
of the general case, or as a strict condition in the definition.
However, some works are not clear in distinguishing between the two situations and this can bring troublesome problems in the future.
In this regard, we would like to mention, for instance, a similar past confusion thoroughly studied in the remarkable paper of
Pach and T\'oth~\cite{DBLP:journals/jct/PachT00} (entitled ``Which crossing number is it anyway?'').

We have demonstrated that a careful distinction (between simple\,/\,non-simple drawings) is surely necessary 
when considering the recent min-$k$-planar graphs.
Our note brings a natural open question about which of the published results of \cite{DBLP:conf/gd/BinucciBBDDHKLMT23}
concerning simply min-$k$-planar graphs with $k\geq2$ remain valid also for general min-$k$-planar graphs.

\bibliography{mink.bib}

\end{document}